\documentclass[journal]{ieeetran}  

\IEEEoverridecommandlockouts                              

\usepackage{graphics} 
\usepackage{epsfig} 
\usepackage{mathptmx} 
\usepackage{times} 
\usepackage{amsmath} 
\usepackage{amssymb}  
\usepackage{cite}
\usepackage{multirow}
\usepackage{color}

\providecommand{\com[2]}{\begin{tt}[#1: #2]\end{tt}}

\providecommand{\rmj[1]}{{\color{black}#1}}

\providecommand{\kcs[1]}{{\color{black}#1}}

\newtheorem{thm}{Theorem}
\newtheorem{problem}{Problem}
\newtheorem{prop}{Proposition}

\newtheorem{rem}{Remark}
\newcommand{\mincutmodif}{{\tt Min Cut with costly nodes} }
\newcommand{\mincut}{{\tt Min Cut} }

\providecommand{\dt}{\Delta \theta}

\title{\LARGE \bf
\kcs{Efficient Computations of a Security Index for False Data Attacks in Power Networks}}

\author{Julien M.~Hendrickx, Karl Henrik Johansson, Raphael M.~Jungers, Henrik Sandberg and Kin Cheong Sou \thanks{\rmj{Julien M. Hendrickx and Raphael M. Jungers are with the Universit\'{e} catholique de Louvain, ICTEAM, Belgium. {\tt \small \{julien.hendrickx,raphael.jungers\}@uclouvain.be}. Their work is supported by the "Communaut\'e francaise de Belgique - Actions de Recherche Concert\'ees", and by the Belgian Programme on Interuniversity Attraction Poles initiated by the Belgian Federal Science Policy Office.  R.M.J. is a F.R.S.-FNRS research associate.} {Kin Cheong Sou, Henrik Sandberg and Karl Henrik Johansson are with the ACCESS Linnaeus Center and the Automatic Control Lab, the School of Electrical Engineering, KTH Royal Institute of Technology, Sweden. {\tt \small \{sou,hsan,kallej\}@kth.se}. Their work is supported by the European Commission through the VIKING project, the Swedish Foundation for Strategic Research (SSF), the Swedish Research Council (VR) and the Knut and Alice Wallenberg Foundation.}}}

\begin{document}

\maketitle
\thispagestyle{empty}
\pagestyle{empty}

%
%
%

\begin{abstract}
The resilience of Supervisory Control and Data Acquisition (SCADA) systems for electric power networks for certain cyber-attacks is considered. We analyze the vulnerability of the measurement system to false data attack on communicated measurements. The vulnerability analysis problem is shown to be NP-hard, meaning that unless $P = NP$ there is no polynomial time algorithm to analyze the vulnerability of the system. Nevertheless, we identify situations, such as the full measurement case, where it can be solved efficiently. In such cases, we show indeed that the problem can be cast as a generalization of the minimum cut problem involving costly nodes. We further show that it can be reformulated as a standard minimum cut problem (without costly nodes) on a modified graph of proportional size. An important consequence of this result is that our approach provides the first exact efficient algorithm for the vulnerability analysis problem under the full measurement assumption. Furthermore, our approach also provides an efficient heuristic algorithm for the general NP-hard problem. Our results are illustrated by numerical studies on benchmark systems including the IEEE 118-bus system.
\end{abstract}


\section{Introduction}
Our society depends heavily on the proper operation of cyber-physical systems, examples of which include, but not limited to, intelligent transport systems, industrial automation systems, health care systems, and electric power distribution and transmission systems. These cyber-physical systems are supervised and controlled through Supervisory Control And Data Acquisition (SCADA) systems. Through remote terminal units (RTUs), SCADA systems collect measurements and send them to the state estimator to estimate the system states. The estimated states are used for subsequent operations such as system health monitoring and control. Any malfunctioning of these operations can lead to significant social and economical consequences such as the northeast US blackout of 2003 \cite{US_blackout2003}.

The technology and the use of the SCADA systems have evolved a lot since they were introduced. The SCADA systems now are interconnected to office LANs, and through them they are connected to the Internet. Hence, today there are more access points to the SCADA systems, and also more functionalities to tamper with. For example, the RTUs can be subjected to denial-of-service attacks. The communicated data can be subjected to false data attacks. Furthermore, the SCADA master itself can be attacked. In the context of secured cyber-physical systems in general, \cite{AminCardenasSastry-HSCC-2009, kn:Gupta2010} have considered denial-of-service-like attacks and their impact. \cite{kn:Bruno09} has studied replay attacks on the sensor measurements and \cite{Smith-IFAC-2011,kn:Pasqualetti2011} have considered false data attacks. This paper investigates the cyber security issues related to false data attacks with the special focus on the measurement systems of power networks. The negative effects of false data attacks on power networks have been exemplified by malware such as Stuxnet and Duqu. False data attacks have received a lot of attention in the literature (e.g.,\ \cite{LRN09,STJ_SCS2010,DS_SGC2010,BRWKNO10,KJTT11,SSJ_CDC2011_mincut,GBMKP11,KP11}). \cite{LRN09} was the first to point out that a coordinated intentional data attack can be staged without being detected by the state estimation bad data detection (BDD) algorithm, a standard part of today's SCADA/EMS system. \cite{LRN09,STJ_SCS2010,DS_SGC2010,KJTT11,SSJ_CDC2011_mincut,GBMKP11,KP11} investigate the construction problem for such ``unobservable'' data attack, especially the sparse ones involving relatively few meters to compromise, under various assumptions of the network (e.g., DC power flow model \cite{Abur_Exposito_SEbook,Monticelli_SEbook}). In particular, \cite{LRN09} poses the attack construction problem as a cardinality minimization problem to find the sparsest attack including a given set of target measurements. References \cite{STJ_SCS2010,DS_SGC2010,SSJ_CDC2011_mincut} set up similar optimization problems for the sparsest attack including a given measurement. References \cite{KJTT11,KP11} seek the sparsest nonzero attack and \cite{GBMKP11} finds the sparsest attack including exactly two injection measurements and possibly more line power flow measurements, under the assumption that all power injections are measured. The solution information of the above optimization problems can help network operators identify the vulnerabilities in the network and strategically assign protection resources (e.g., encryption of measurements and secured PMUs) to their best effect (e.g., \cite{DS_SGC2010,BRWKNO10,KP11}). On the other hand, the unobservable data attack problem has its connection to another vital EMS functionality -- observability analysis \cite{Abur_Exposito_SEbook,Monticelli_SEbook}. In particular, solving the attack construction problem can also solve an observability analysis problem (this was explained in \cite[Section II-C]{SSJ_exact_l1}). This connection was first reported in \cite{KJTT11}, and was utilized in \cite{SSJ_ckt} to compute the sparsest critical $p$-tuples for some integer $p$. This is a generalization of critical measurements and critical sets \cite{Abur_Exposito_SEbook}.

To perform the cyber-security analysis in a timely manner, it is important to solve the data attack construction problem efficiently. This effort has been discussed, for instance, in \cite{LRN09,STJ_SCS2010,DS_SGC2010,KJTT11,SSJ_CDC2011_mincut,GBMKP11,KP11,SSJ_exact_l1}. The efficient solution to the attack construction problem in \cite{STJ_SCS2010} is the focus of this paper. The matching pursuit method \cite{Mallat93matchingpursuit} employed in \cite{LRN09} and the basis pursuit method \cite{CDS_Basis_Pursuit} ($l_1$ relaxation and its weighted variant) employed in \cite{KP11} are common efficient (i.e., polynomial time) approaches to suboptimally solve the attack construction problem. However, these methods do not guarantee exact optimal solutions, and they might not be sufficiently accurate. For instance, \cite{SSJ_CDC2011_mincut} describes a naive application of basis pursuit and its consequences. While \cite{KJTT11,KP11} provide polynomial time solution procedures for their respective attack construction problems, the problems therein are different from the one in this paper in that the considered problem in this paper is not a special case of the ones in \cite{KJTT11,KP11}. In particular, in \cite{KJTT11} the attack vector contains at least one nonzero entry. However, this nonzero entry cannot be given a priori. This means that the problem considered in this paper is more general than the one in \cite{KJTT11}. \cite{KP11} needs to restrict the number of nonzero injection measurements attacked, while there is no such constraint in the problem considered in this paper.

Other relevant previous work include \cite{SSJ_ckt,SSJ_CDC2011_mincut,SSJ_exact_l1}, which also consider the data attack construction problem in this paper. In \cite{SSJ_ckt,SSJ_CDC2011_mincut} the attack construction problem is formulated as a graph generalized minimum cut problem (to be defined in Section~\ref{subsec:gmc_costly_nodes}). However, it is not known in \cite{SSJ_ckt,SSJ_CDC2011_mincut} whether the generalized minimum cut problem can be solved efficiently (i.e., in polynomial time) or not. Indeed, \cite{SSJ_ckt,SSJ_CDC2011_mincut} only provide approximate solutions. Instead, the current work establishes that the generalized minimum cut problem is indeed \emph{exactly solvable in polynomial time}. This work establishes the result by constructing a practical polynomial time algorithm. Regarding \cite{SSJ_exact_l1}, one of the main distinctions is that \cite{SSJ_exact_l1} makes an assumption that no bus injections are metered. The current result requires a different assumption that the network is \emph{fully measured} as in \cite{GBMKP11} (i.e., all bus injections and line power flows are metered). In addition, \cite{SSJ_exact_l1} considers a more general case where the constraint matrix is totally unimodular, whereas the focus of the current paper is a graph problem. The setup considered by this paper is specific to power network applications and thus it enables a more efficient solution algorithm. Finally, note that the notion of minimum cut problem has been explored also in other power network applications (e.g., \cite{LRDP06,inhibitbisect,ETEP:ETEP255}).

{\bf Outline:} In the next section, we present the optimization problem of interest, namely the security index problem, and discuss its applications. Then Sections~\ref{sec:complexity}, \ref{sec:si_full_measurement} and \ref{section-reformulation} present the technical contributions of this paper, focusing on a specialized version of the security index problem defined in (\ref{opt:security_index}) in the end of Section~\ref{subsec:security_index}. In Section~\ref{sec:complexity} the complexity of the security index problem is analyzed. We show that the security index problem is NP-hard in general, but in Section~\ref{sec:si_full_measurement} we demonstrate that under some realistic assumptions it can be restated as a generalized minimum cut ({\tt Min Cut}) problem. In Section~\ref{section-reformulation} we show that the generalized \mincut problem can be solved efficiently, by reformulating it as a classical \mincut problem. The specialized version considered in Sections~\ref{sec:complexity}, \ref{sec:si_full_measurement} and \ref{section-reformulation} turns out to be not restrictive, as far as the application of the proposed results is concerned. This will be explained in Section~\ref{sec:si_edge_only}. In Section~\ref{section-example} a simple numerical example is first presented to illustrate that the proposed solution correctly solves the generalized \mincut problem, while previous methods cannot. Then the efficiency and accuracy of the proposed solution to the security index problem are demonstrated through a case study with large-scale benchmark systems. We also demonstrate that our method provides an efficient and high quality approximate solution to the general problem security index problem which is NP-hard. 


\section{The Security Index Problem}
In Section~\ref{subsec:model} the mathematical model of the power networks considered is first described. Then in Section~\ref{subsec:security_index} the security index of power networks is defined.

\subsection{Power Network Model and State Estimation} \label{subsec:model}
A power network is modeled as a graph with $n+1$ nodes and $m$ edges. The nodes and edges model the buses and transmission lines in the power network, respectively. In the present text, the terms node and bus are used interchangeably, and the same is true for edges and transmission lines (or simply lines). The topology of the graph is described by a directed incidence matrix $A \in \mathbb{R}^{(n+1) \times m}$, in which the directions along the edges are arbitrarily specified \cite{SSJ_CDC2011_mincut}. The physical property of the network is described by a nonsingular diagonal matrix $D \in \mathbb{R}^{m \times m}$, whose nonzero entries are the reciprocals of the reactance of the transmission lines. In general, the reactance is positive (i.e., inductive) and hence the matrix $D$ is assumed to be positive definite throughout this paper. 

In the sequel, the set of all nodes and the set of all directed edges of the power network graph are denoted $V^0$ and $E^0$, respectively. The edge directions are consistent with those in $A$. An element of $V^0$ is denoted by $v_i \in V^0$, and an element of $E^0$ is denoted by $(v_i, v_j) \in E^0$ for $v_i \in V^0$ and $v_j \in V^0$. The set of all neighbors of $v_i$ is denoted by $N(v_i)$. A node $v_j$ is a neighbor of $v_i$ if either $(v_i, v_j) \in E^0$ or $(v_j, v_i) \in E^0$.

The states of the network include bus voltage phase angles and bus voltage magnitudes, the latter of which are typically assumed to be constant (one in the per unit system). Therefore, the network states can be captured in a vector $\theta \in {[0,2\pi)}^{n+1}$. The state estimator estimates $\theta$ based on the measurements obtained from the network. In reality the model relating the states and the measurements is nonlinear. However, for state estimation data attack analysis \cite{LRN09,STJ_SCS2010,DS_SGC2010,BRWKNO10,KJTT11,SSJ_CDC2011_mincut,GBMKP11,KP11,SSJ_exact_l1} (and more traditionally bad data analysis \cite{Abur_Exposito_SEbook,Monticelli_SEbook,HSKF75}) it suffices to consider the DC power flow model \cite{Abur_Exposito_SEbook,Monticelli_SEbook}. In the DC power flow model the measurement vector, denoted as $z$, is related to $\theta$ by
\begin{equation} \label{def:H_matrix}
  z = H \theta + \Delta z, \quad {\rm where} \quad H \triangleq \begin{bmatrix}
    P_1 D A^T \\ -P_2 D A^T \\ P_3 A D A^T
  \end{bmatrix}.
\end{equation}
In (\ref{def:H_matrix}), $\Delta z$ can either be a vector of random error or intentional additive data attack (e.g., \cite{LRN09}), and $P_1$, $P_2$ and $P_3$ consist of subsets of rows of identity matrices of appropriate dimensions, indicating which measurements are actually taken. The term $P_1 D A^T \theta$ contains line power flow measurements, measured at the outgoing ends of the lines. Similarly, $-P_2 D A^T \theta$ contains the line power flow measurements at the incoming ends of the lines. The term $P_3 A D A^T \theta$ contains bus power injection measurements, one entry for each measured bus.

Measurement redundancy is a common practice in power networks \cite{Abur_Exposito_SEbook,Monticelli_SEbook}. Therefore, it is assumed in this paper that the measurement system described by $H$ is observable -- meaning that if any column of $H$ is removed the remaining submatrix still has rank $n$ \cite{Abur_Exposito_SEbook,Monticelli_SEbook}. Note that $H$ cannot have rank $n+1$ since the sum of all columns of $H$ is always a zero column vector (a property of any incidence matrix $A$). In the practice of power system state estimation, it is customary to designate an arbitrary node as the reference and set the corresponding entry of $\theta$ to zero. Without loss of generality, it is assumed that the first entry of $\theta$ is zero (i.e., $\theta(1) = 0$) and denote $\theta_{2:}$ as the rest of the entries of $\theta$. For convenience, let $H_{2:}$ denote $H$ with the first column removed. By definition, $H \theta = H_{2:} \theta_{2:}$ and $H_{2:}$ has full column rank ($=n$) since $H$ is observable. Given measurements $z$, the estimate of the network states is typically determined via the least squares approach \cite{Abur_Exposito_SEbook,Monticelli_SEbook}:
\begin{equation} \label{def:theta_hat}
  \hat{\theta}_{2:} = {({H_{2:}}^T W H_{2:})}^{-1} {H_{2:}}^T W z,
\end{equation}
where $W$ is a given positive-definite diagonal matrix, whose nonzero entries are typically the reciprocals of the variances of the measurement noise. The state estimate $\hat{\theta} = {\begin{bmatrix} 0 & {\hat{\theta}_{2:}}^T \end{bmatrix}}^T$ is subsequently fed to other vital SCADA functionalities such as optimal power flow (OPF) calculation and contingency analysis (CA). Therefore, the accuracy and reliability of $\hat{\theta}$ is of paramount concern.

\subsection{Security Index} \label{subsec:security_index}
To detect possible faults or data attacks in the measurements $z$, the BDD test is commonly performed (\cite{Abur_Exposito_SEbook,Monticelli_SEbook}). In a typical strategy, if the norm of the residual
\begin{equation} \label{def:BDD_residual_norm}
  {\rm residual} \triangleq z - H_{2:} \hat{\theta}_2 = (I - H_{2:}{({H_{2:}}^T W H_{2:})}^{-1} {H_{2:}}^T W) \Delta z
\end{equation}
is too big, then the BDD alarm is triggered. The BDD test is in general sufficient to detect the presence of a random error $\Delta z$ \cite{Abur_Exposito_SEbook,Monticelli_SEbook}. However, in face of a coordinated malicious attack the BDD test can fail. In particular, in \cite{LRN09} it was reported that an attack of the form
\begin{equation} \label{def:unobservable_attack}
\Delta z = H \Delta \theta
\end{equation}
for an arbitrary $\Delta \theta \in \mathbb{R}^{n+1}$ would result in a zero residual in (\ref{def:BDD_residual_norm}) since $H \Delta \theta = H_{2:} \Delta \theta_{2:}$ for some $\Delta \theta_{2:} \in \mathbb{R}^n$. Data attack in the form of (\ref{def:unobservable_attack}) is unobservable from the BDD perspective, and this was also experimentally verified in \cite{TDSJ11} in a realistic SCADA system testbed. Since \cite{LRN09}, there has been a significant amount of literature studying the unobservable attack in (\ref{def:unobservable_attack}) and its consequences to state estimation data integrity (e.g., \cite{STJ_SCS2010,DS_SGC2010,BRWKNO10,KJTT11,GBMKP11,KP11}). In particular, \cite{STJ_SCS2010} introduced the notion of security index $\alpha_k$ for a measurement $k$ as the optimal objective value of the following cardinality minimization problem:
\begin{equation} \label{opt:card_min_con}
  \begin{array}{cccl}
    \alpha_k & \triangleq & \mathop{\rm min}\limits_{\Delta \theta \in \mathbb{R}^{n+1}} & {\rm card}\big( H \Delta \theta \big) \vspace{1mm} \\
    & & \textrm{subject to} & H(k,:) \Delta \theta = 1,
  \end{array}
\end{equation}
where ${\rm card}(\cdot)$ denotes the cardinality of its argument, $k$ is the label of the measurement for which the security index $\alpha_k$ is computed, and $H(k,:)$ denotes the $k^{\rm th}$ row of $H$. $\alpha_k$ is the minimum number of measurements an attacker needs to compromise in order to attack measurement $k$ without being detected. In particular, a small $\alpha_k$ implies that measurement $k$ is relatively easy to compromise in an unobservable attack. As a result, the knowledge of the security indices for all measurements allows the network operator to pinpoint the security vulnerabilities of the network, and to better protect the network with limited resource. For example, \cite{DS_SGC2010} proposed a method to optimally assign limited encryption protection resources to improve the security of the network based on its security indices.

It should be emphasized that the security index defined in (\ref{opt:card_min_con}) can provide a security assessment that the standard power network BDD procedure \cite{Abur_Exposito_SEbook,Monticelli_SEbook} might not be able to provide. As a concrete example \cite{STJ_SCS2010}, consider the simple network whose $H_{2:}$ matrix is
\begin{equation}
H_{2:} = \begin{pmatrix} -1 & -1 & 0 \\ -1 & 0 & 0 \\ 1 & 0 & 0 \\
1 & 0 & -1 \\ 0 & -1 & 0
\end{pmatrix}.
\label{eq:Hex}
\end{equation}
From (\ref{def:theta_hat}), the ``hat matrix'' \cite{Abur_Exposito_SEbook,Monticelli_SEbook}, denoted $K$ is defined according to
\begin{displaymath}
  \hat{z} = H_{2:} \hat{\theta}_2 = H_{2:} {({H_{2:}}^T W H_{2:})}^{-1} {H_{2:}}^T W z \triangleq K z.
\end{displaymath}
Assuming $W = I$, the $K$ matrix associated with $H_{2:}$ in (\ref{eq:Hex}) is
\begin{equation}
K = \begin{pmatrix}  0.6&    0.2 &   -0.2&         0   & 0.4 \\
    0.2 &    0.4 &   -0.4 &         0   & -0.2 \\
   -0.2 &   -0.4 &    0.4 &         0   &  0.2 \\
         0 &        0 &         0 &    1 &         0 \\
    0.4 &   -0.2 &    0.2 &         0 &    0.6
\end{pmatrix}.
\label{eq:Kex}
\end{equation}
The hat matrix $K$ shows how the measurements $z$ are
weighted together to form a power flow estimate $\hat{z}$. The rows
of the hat matrix can be used to study the measurement redundancy
in the system \cite{Abur_Exposito_SEbook,Monticelli_SEbook}. Typically a large degree of
redundancy (many non-zero entries in each row) is desirable to
compensate for noisy or missing measurements. In (\ref{eq:Kex}),
it is seen that all measurements are redundant except the
fourth which is called a \emph{critical
measurement}. Without the critical measurement observability is
lost. From the hat matrix one is led to believe that the critical
measurement is sensitive to attacks. This is indeed the case, but some other measurements can also be vulnerable to attacks. It can be shown - for example using the method that we develop, that the security indices $\alpha_k$, $k = 1,\ldots,5$, respectively, are 2, 3, 3, 1, 2. Therefore, the fourth measurement (critical measurement) has security index one, indicating that it is indeed vulnerable to unobservable attacks. However, the first and the last measurements also have relatively small security indices. This is not obvious from $K$ in (\ref{eq:Kex}). Hence, we cannot rely on the hat matrix for vulnerability analysis of power networks.

For ease of exposition but without loss of generality, instead of (\ref{opt:card_min_con}) the following version of the security index problem with a specialized constraint will be the focus of the parts of the paper where the main technical contributions are presented (i.e., Sections~\ref{sec:complexity}, \ref{sec:si_full_measurement} and \ref{section-reformulation}):
\begin{equation} \label{opt:security_index}
  \begin{array}{cl}
    \mathop{\rm minimize}\limits_{\Delta \theta \in {\mathbb{R}}^{n+1}} & {\rm card}\big( H \Delta \theta \big) \vspace{1mm} \\
    \textrm{subject to} & A(:,e)^T \Delta \theta \neq 0,
  \end{array}
\end{equation}
where $e \in \{1,2,\ldots,m\}$ is given. The restriction introduced in (\ref{opt:security_index}) is that it can only enforce constraints on edge flows but not on node injections as directly allowed by (\ref{opt:card_min_con}). We will see however in Section~\ref{sec:si_edge_only} that all results obtained for (\ref{opt:security_index}) can be extended to the general case in (\ref{opt:card_min_con}).



\section{The Security Index Problem is NP-hard} \label{sec:complexity}
Consider a variant of (\ref{opt:card_min_con}) where $k$ is not fixed (i.e., one wishes to minimize ${\rm card} \left(H \Delta \theta\right)$ under the constraint that at least one entry of $H \Delta \theta$ is nonzero). This variant of (\ref{opt:card_min_con}) is known to be the cospark of $H_{2:}$ in compressed sensing \cite{CT05}. The cospark of $H_{2:}$ is the same as the spark of $F$, where $F$ is a matrix of full row rank such that $F H_{2:} = 0$ \cite{CT05}. The spark of $F$ is defined as the minimum number of columns of $F$ which are linearly dependent \cite{Donoho04032003}. It is established that computing the spark of a general matrix $F$ is NP-hard \cite{spark_NP_hard,McCormick_PhD}. Consequently, because of the equivalence between spark and cospark, unless $P=NP$ there is no efficient algorithm to solve the security index problem in (\ref{opt:card_min_con}) if the $H$ matrix is not assumed to retain any special structure. In power network applications, the $H$ matrix in fact possesses special structure as defined in (\ref{def:H_matrix}). Nevertheless, the security index problem, even the specialized version in (\ref{opt:security_index}), is still computationally intractable as indicated by the following statement:
\begin{thm}\label{thm:NP_hardness}
Unless $P=NP$, there is no polynomial time algorithm that solves the problem \eqref{opt:security_index}, with $H$ defined in \eqref{def:H_matrix}, even if $D$ is the identity matrix and $P_2=0$.
\end{thm}
\begin{IEEEproof}
  Our proof proceeds by reduction from the positive one-in-three 3SAT problem \cite{gopalan2006connectivity}:
\emph{Given a set of $M$ triples of indices $C_j = (\alpha_j,\beta_j,\gamma_j)\in \{1,\dots,n\}^3$, does there exist a vector $\tilde x\in \{0,1\}^n$ such that for every $j$, exactly one among $\tilde x_{\alpha_j}, \tilde x_{\beta_j}, \tilde x_{\gamma_j}$ is 1 and the others 0.}

Consider an instance of the positive one-in-three 3SAT problem, and let us build an equivalent instance of (\ref{opt:card_min_con}). We set $P_2$ to 0, and set $D$ as the identity matrix. As a result, non-trivially zero entries of $H\dt$ corresponding to edges $(i,j)$ will be of the form $\dt_i - \dt_j$, while those corresponding to a node $i$ will be of the form $\sum_{j:(i,j)\in E} (\dt_j-\dt_i)$. We remind that the entry of an edge is trivially zero if the corresponding entry in $P_1$ is 0, and that of a node is trivially zero if the corresponding entry of $P_3$ is 0.

We begin by defining a node 1 and a node 0 connected by an edge whose corresponding entry in $P_1$ is set to 1. We set $k$ such that the constraint $H(k,:)\dt = 1$ in \eqref{opt:card_min_con}  corresponds to this edge, so that their must hold $\dt_1-\dt_0 = 1$ for any solution of the problem. Since $H\dt$ is not modified when adding a constant to all entries of $\dt$, we assume without loss of generality that $\dt_1 =1$ and $\dt_0 =0$.

The goal of the  first part of our construction is to represent the variables. For every $i=1,\dots,n$, we define a node $x_i$ that we connect to both 1 and 0. We set to 1 the entries of $P_1$ corresponding to the edges $(1,x_i)$ and $(0,x_i)$, and to 0 the entries of $P_3$ corresponding to $x_i$.
Observe that the two entries of $H\dt$ corresponding to these two edges are $1-\dt_{x_i}$ and $0-\dt_{x_i}$, which cannot be simultaneously 0. Moreover, one of them is equal to zero if and only if $\dt_{x_i}$ is either $0$ or $1$.

Taking into account the fact the entry of $H\dt$ corresponding to the edge $(1,0)$ is by definition 1, we have thus proved that ${\rm card}\left(H\dt\right) \geq n+1$ for any $\dt$, independently of the rest of the construction. Moreover, ${\rm card}\left(H\dt\right)=n+1$  only if  $\dt_{x_i}\in \{0,1\}$ for every $i$, and if the entries of $H\dt$ corresponding to all the edges and nodes introduced in the sequel are 0.
The remainder of the construction, represented in Fig.~\ref{fig:construction_NP}, is designed to ensure that all these entries can be 0 only if the (binary) values $\dt_{x_i}$ solves the initial instance of the positive one-in-three 3SAT problem.

We first generate a reference value at $\frac{1}{3}$ for every clause: We define two nodes indexed by $\frac{2}{3}$ and $\frac{1}{3}$, and add the connections $(1,\frac{2}{3}),(\frac{2}{3},\frac{1}{3}),(\frac{1}{3},0)$. The entries of $P_1$ corresponding to these connections are set to 0, but the entry of $P_3$ corresponding to the nodes $\frac{1}{3}$ and $\frac{2}{3}$ are set to 1.
Besides, we define for every clause $j=1,\dots,M$ a clause node $c_j$ connected to $\frac{1}{3}$ by an edge whose corresponding entry in $P_1$ is set to one.

The entries of $H\dt$ corresponding to the edges between $\frac{1}{3}$ and $c_j$  are $\dt_{\frac{1}{3}}-\dt_{c_j}$, which are thus zero only when
$\dt_{c_j}=\dt_{\frac{1}{3}}$ for every $j$. Using these equalities, observe now that the entry of $H\dt$ corresponding to $\frac{1}{3}$ is
$$
\dt_{\frac{2}{3}} + \dt_{0} + \sum_{j=1}^M \dt_{c_j} - (2+M) \dt_{\frac{1}{3}} = \dt_{\frac{2}{3}} + 0 - 2 \dt_\frac{1}{3},
$$
while the entry corresponding to $\frac{2}{3}$ is $1+\dt_{\frac{1}{3}} - 2\dt_{\frac{2}{3}}$. These two entries are thus equal to zero if and only if $\dt_{\frac{2}{3}} = \frac{2}{3}$ and $\dt_{c_j} =\dt_{\frac{1}{3}}= \frac{1}{3}$ for every $j$, as intended.

We now represent the clauses. 
For each $j$, we connect the clause node $c_j$ to the nodes $x_{\alpha_j}$, $x_{\beta_j}$ and $x_{\gamma_j}$ of the three variables involved by edges whose corresponding entries in $P_1$ are zero. On the other hand, we set to 1 the entry of $P_3$ corresponding to $c_j$. The corresponding (non trivially zero) entries of $H\dt$ are then
$$
\dt{x_{\alpha_j}}+ \dt{x_{\beta_j}}+ \dt{x_{\gamma_j}}-3\dt_{c_j} =
\dt{x_{\alpha_j}}+ \dt{x_{\beta_j}}+ \dt{x_{\gamma_j}}-1.
$$
Remembering that $\dt_{x_i}$ is either 1 or 0 for any $i$, this latter expression can be zero only if exactly one among $\dt{x_{\alpha_j}}$, $\dt{x_{\beta_j}}$ and $\dt{x_{\gamma_j}}$ is 1.
If that is the case, setting $\tilde x_i= \dt_{x_i}$ for every $i$ yields a vector $x$ that solves the instance of positive one-in-three 3SAT.

We have thus shown that there exists a $\dt$ for which ${\rm card}\left(H \dt\right)=n+1$ only if the $\dt_{x_i}$ are binary, and if the binary vector $\tilde x$ obtained by setting $\tilde x_i = \dt_{x_i}$ solves the instance of positive one-in-three 3SAT. Conversely, one can verify that if a binary vector $\tilde x$ solves the instance of the one-in-three 3SAT problem, then setting $\dt_{x_i} = \tilde x_i$ for every $i$, $\dt_{c_j}=\dt_{\frac{1}{3}} = \frac{1}{3}$ for every  $j$  and $\dt_{\frac{2}{3}} = \frac{2}{3}$ yields a cost ${\rm card}\left(H\dt\right)=n+1$. The latter cost can thus be obtained if and only if the initial positive one-in-three 3SAT problem is achievable. This achieves the proof because our construction clearly takes an amount of time that grows polynomially with the size of the instance $C$, and unless $P=NP$ there is no polynomial time algorithm that solves the positive one-in-three 3SAT \cite{gopalan2006connectivity}.

\begin{figure}\centering
\includegraphics[scale=.35]{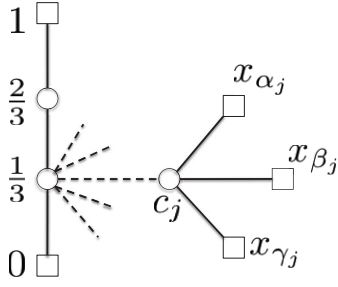}
\caption{Representation of a part of the construction of the proof of Theorem \ref{thm:NP_hardness}, including the reference values of $\dt$ and one clause $C_j$. Edges are represented by dashed line when they are measured and continuous lines otherwise. Nodes are represented by squares when they are measured and circles otherwise.
If ${\rm card}\left(H\dt\right)=n+1$, $\dt$ takes only values 1 and 0 for the $x_i$ and all  entries of $H\dt$ other than those corresponding to the nodes $x_i$ must be zero. As a result, a dashed edge transmits no current and enforces equality between the values of the nodes to which it is incident, and circle nodes enforce that the sum of the currents on the incident edges should be 0. These constraints can only be satisfied if $\dt_{c_j}=\frac{1}{3}$, and if exactly one of the nodes involved in each clause is at 1 and the others at 0.}\label{fig:construction_NP}
\end{figure}
\end{IEEEproof}
\begin{rem}
  (\ref{opt:card_min_con}) is also NP-hard since (\ref{opt:security_index}) is a special case of (\ref{opt:card_min_con}).
\end{rem}

\section{Tractable Special Cases of the Security Index Problem} \label{sec:si_full_measurement}

In Section~\ref{subsec:full_measurement} we show that, under the full measurement assumption, the security index problem can be solved by solving its restriction where decision variables take binary values. Section~\ref{subsec:equiv_general_thm} presents the proof of the statement which implies our finding in Section~\ref{subsec:full_measurement}. Section~\ref{subsec:equiv_general_thm} also discusses the relationship between the security index problem and its binary restriction defined in Section~\ref{subsec:full_measurement}. Section~\ref{subsec:gmc_costly_nodes} describes the consequences of Sections~\ref{subsec:full_measurement} and \ref{subsec:equiv_general_thm}, explaining how the security index problem can be reformulated as a generalized minimum cut problem with costly nodes, a graph problem whose efficient solution will be discussed in Section~\ref{section-reformulation}.

\subsection{The Security Index Problem Under Full Measurement Assumption} \label{subsec:full_measurement}
Even though in general the security index problem in (\ref{opt:security_index}) is NP-hard for $H$ defined in (\ref{def:H_matrix}), there exist interesting specializations that are solvable in polynomial time. One such case is the \emph{full measurement} situation where $P_1 = I$, $P_2 = I$ and $P_3 = I$. In \cite{SSJ_CDC2011_mincut,GBMKP11} the full measurement assumption is also considered, motivated by the situations where all power flows and injections are measured in future smart grid applications. The polynomial time complexity of (\ref{opt:security_index}) under the full measurement assumption can be established in three steps. Firstly, it can be shown that problem~(\ref{opt:security_index}) can be solved by solving a restriction where the decision vector $\Delta \theta$ is a binary vector. Secondly, in Section~\ref{subsec:gmc_costly_nodes} it will be shown that the binary restriction of (\ref{opt:security_index}) can be expressed in a generalized \mincut problem with costly nodes. Finally, this generalized \mincut problem can be shown to be solvable in polynomial time. This is to be explained in Section~\ref{section-reformulation}. 

The first step is formalized in the following statement, whose preliminary version appeared in \cite{SSJ_CDC2011_mincut}.
\begin{prop} \label{thm:si_sibin_equiv}
  Let $H$ in (\ref{def:H_matrix}) satisfy the full measurement assumption that $P_1 = I$, $P_2 = I$, and $P_3 = I$. Consider the following restriction of problem~(\ref{opt:security_index}) with 0-1 binary decision vector:
  \begin{equation} \label{opt:security_index_binary}
  \begin{array}{cl}
    \mathop{\rm minimize}\limits_{\Delta \theta \in {\left\{0,1\right\}}^{n+1}} & {\rm card}\big( H \Delta \theta \big) \vspace{1mm} \\
    \textrm{subject to} & {A(:,e)}^T \Delta \theta \neq 0.
  \end{array}
\end{equation}
  It holds that every optimal solution of (\ref{opt:security_index_binary}) is an optimal solution of (\ref{opt:security_index}) (i.e., the problem with the same formulation except that  $\Delta \theta$ is not restricted to binary values).
\end{prop}
\begin{IEEEproof}
  Proposition~\ref{thm:si_sibin_equiv} is a corollary of the more general Theorem~\ref{thm:si_sibin_equiv_general} to be described in Section~\ref{subsec:equiv_general_thm}.
\end{IEEEproof}
\begin{rem}
  Since there cannot be any all zero column in any incidence matrix $A$, problem~({\ref{opt:security_index}}) and (\ref{opt:security_index_binary}) are always feasible. Proposition~\ref{thm:si_sibin_equiv} states that, under the full measurement assumption, an optimal solution of ({\ref{opt:security_index}}) can always be obtained by solving (\ref{opt:security_index_binary}). The later problem will be shown to be solvable in polynomial time.
\end{rem}

\subsection{The Security Index Problem with Binary Decision Vector}
\label{subsec:equiv_general_thm}

In the sequel, let $c_{ij} \ge 0$ represent the cost of attacking the power flow measurements of a line $(v_i,v_j)$, and $p_i \ge 0$ the cost of attacking the injection measurement at bus $v_i$. Problems (\ref{opt:security_index}) can be reformulated in a more general way that also allows taking into account the fact that tempering with certain measurements may be more expensive than with some others:
\begin{equation} \label{opt:gsi}
    \begin{array}{cl}
      \mathop{\rm minimize}\limits_{\Delta \theta \in {\mathbb{R}}^{n+1}} & c^T g(D A^T \Delta \theta) + p^T g(A D A^T \Delta \theta) \vspace{1mm} \\
      \textrm{subject to} & {A(:,e)}^T \Delta \theta \neq 0,
    \end{array}
\end{equation}
and the corresponding reformulation of (\ref{opt:security_index_binary}) is defined by
\begin{equation} \label{opt:gsi_bin}
    \begin{array}{cl}
      \mathop{\rm minimize}\limits_{\Delta \theta \in {\left\{0,1\right\}}^{n+1}} & c^T g(D A^T \Delta \theta) + p^T g(A D A^T \Delta \theta) \vspace{1mm} \\
      \textrm{subject to} & {A(:,e)}^T \Delta \theta \neq 0.
    \end{array}
\end{equation}
In (\ref{opt:gsi}) and (\ref{opt:gsi_bin}), $g$ is a vector-valued indicator function such that for any vector $x$, $g_i(x) = 1$ if $x_i \neq 0$ and $g_i(x) = 0$ otherwise. It can be seen that if $c_{ij} \in \{0,1,2\}$ and $p_i \in \{0,1\}$, then (\ref{opt:security_index}) and (\ref{opt:security_index_binary}) are recovered. Let ${\begin{bmatrix} D A^T \Delta \theta \end{bmatrix}}_{(v_i,v_j)}$ denote the entry of $D A^T \Delta \theta$ corresponding to edge $(v_i,v_j)$, and let ${\begin{bmatrix} A D A^T \Delta \theta \end{bmatrix}}_{v_i}$ denote the entry of $A D A^T \Delta \theta$ corresponding to node $v_i$. With a slight abuse of notation, the symbol $g\big({\begin{bmatrix} D A^T \Delta \theta\end{bmatrix}}_{(v_i,v_j)} \big)$ denotes the entry of $g(D A^T \Delta \theta)$ corresponding to $(v_i,v_j)$. In addition, $g\left({\begin{bmatrix} A D A^T \Delta \theta \end{bmatrix}}_{v_i} \right)$ is defined similarly.

The following theorem characterizes the relationship between the security index problem in (\ref{opt:security_index}) and its binary restriction in (\ref{opt:security_index_binary}) by studying their respective generalizations of (\ref{opt:gsi}) and (\ref{opt:gsi_bin}) for arbitrary nonnegative vectors $c$ and $p$.
\begin{thm} \label{thm:si_sibin_equiv_general}
  Let $J_c$ and $J_b$, respectively, denote the optimal objective values of (\ref{opt:gsi}) and (\ref{opt:gsi_bin}) with $A$ and $D$ defined in (\ref{def:H_matrix}), $c \in \mathbb{R}^m_+$, $p \in \mathbb{R}^{n+1}_+$ and $e \in \{1,2,\ldots,m\}$ given. Then
  \begin{equation} \label{eqn:si_sibin_error_bound}
    0 \le J_b - J_c \le \sum\limits_{v_i \in V^0} \max\Big\{0, \max\limits_{v_j \in N(v_i)} \big\{p_i - c_{ij} \big\} \Big\}.
  \end{equation}
\end{thm}
\begin{IEEEproof}
  First note that both (\ref{opt:gsi}) and (\ref{opt:gsi_bin}) are always feasible with finite optimal objective values attained by some optimal solutions. In addition, $0 \le J_b - J_c$ holds because (\ref{opt:gsi_bin}) is a restriction of (\ref{opt:gsi}). To show the upper bound in (\ref{eqn:si_sibin_error_bound}) the main idea is that for each feasible solution $\Delta \theta$ of (\ref{opt:gsi}) it is possible to construct a feasible solution $\Delta \phi$ of (\ref{opt:gsi_bin}), such that the objective value difference is bounded from above by $\sum\limits_{v_i \in V^0} \max\Big\{0, \max\limits_{v_j \in N(v_i)} \big\{p_i - c_{ij} \big\} \Big\}$. The construction is as follows. Let $\Delta \theta$ be a feasible solution of (\ref{opt:gsi}), and let $\Delta \theta(v_i)$ be its entry corresponding to node $v_i \in V^0$. Since $\Delta \theta$ is feasible, the constraint $A(:,e)^T \Delta \theta \neq 0$ implies that there exist two nodes denoted $v_s$ and $v_t$ with $e$ corresponding to either $(v_s,v_t)$ or $(v_t,v_s)$ such that $\Delta \theta(v_s) \neq \Delta \theta(v_t)$. Without loss of generality, it is assumed that $\Delta \theta(v_s) > \Delta \theta(v_t)$. Define $\Delta \phi \in {\{0,1\}}^{n+1}$ by
  \begin{equation} \label{def:delta_tilde_theta}
\Delta \phi(v_i) = \begin{cases}
  1 & \text{if} \; \Delta \theta(v_i) > \Delta \theta(v_t) \\
  0 & \text{if} \; \Delta \theta(v_i) \le \Delta \theta(v_t)
\end{cases}
\quad \forall \; v_i \in V^0.
\end{equation}
  Note that $\Delta \phi$ is feasible to (\ref{opt:gsi_bin}) since $\Delta \phi(v_s) \neq \Delta \phi(v_t)$ by construction. Also notice that for any two nodes $v_i$ and $v_j$ if $\Delta \theta(v_i) = \Delta \theta(v_j)$ then $\Delta \phi(v_i) = \Delta \phi(v_j)$. Hence, in the objective functions of (\ref{opt:gsi}) and (\ref{opt:gsi_bin}) it holds that
  \begin{equation} \label{eqn:pf_edge_cost_noincrease}
    c_{ij} g\big({\begin{bmatrix} D A^T \Delta \theta\end{bmatrix}}_{(v_i,v_j)} \big) \ge c_{ij} g\big({\begin{bmatrix} D A^T \Delta \phi\end{bmatrix}}_{(v_i,v_j)} \big), \quad \forall \; (v_i,v_j) \in E^0.
  \end{equation}
  In other words, for each edge the contribution to the objective function with the new solution $\Delta \phi$ is smaller than or equal to that with the initial one $\Delta \theta$. To finish the proof, the objective function contribution due to the node injections needs to be investigated. Let $V_b \subset V^0$ be defined such that
  \begin{equation} \label{eqn:pf_Vb_def}
    v_i \in V_b \;\; \iff \;\; g\left({\begin{bmatrix} A D A^T \Delta \theta \end{bmatrix}}_{v_i} \right) = 0, \; \; g\left({\begin{bmatrix} A D A^T \Delta \phi \end{bmatrix}}_{v_i}\right) = 1.
  \end{equation}
  In essence, $V_b$ encompasses all causes for $J_b > J_c$. Consider $v_i \in V_b$, since $g\left({\begin{bmatrix} A D A^T \Delta \phi \end{bmatrix}}_{v_i}\right) = 1$, there exists $v_k \in N(v_i)$ such that $\Delta \theta(v_k) \neq \Delta \theta(v_i)$. Consequently, the fact that $g \left({\begin{bmatrix} A D A^T \Delta \theta \end{bmatrix}}_{v_i}\right) = 0$ implies that there exists $v_i^+ = {\rm argmax}_{v_k \in N(v_i)} \{\Delta \theta(v_k)\}$ such that $\Delta \theta(v_i^+) > \Delta \theta(v_i)$. Similarly, there exists $v_i^- = {\rm argmin}_{v_k \in N(v_i)} \{\Delta \theta(v_k)\}$ such that $\Delta \theta(v_i^-) < \Delta \theta(v_i)$. If $\Delta \theta(v_i) > \Delta \theta(v_t)$, then (\ref{def:delta_tilde_theta}) implies that $\Delta \phi(v_i^+) = \Delta \phi(v_i) = 1$. In addition, it holds that if $v_i^+ \in V_b$ then $v_i \neq {\rm argmax}_{v_k \in N(v_i^+)} \{\Delta \theta(v_k)\}$. This is true because ${\rm argmax}_{v_k \in N(v_i^+)} \{\Delta \theta(v_k)\} > \Delta \theta(v_i^+) > \Delta \theta(v_i)$ if $v_i^+ \in V_b$. Conversely, if $\Delta \theta(v_i) \le \Delta \theta(v_t)$, then $\Delta \phi(v_i^-) = \Delta \phi(v_i) = 0$. Similar to the case with $v_i^+$, if $v_i^- \in V_b$, then $v_i \neq {\rm argmin}_{v_k \in N(v_i^-)} \{\Delta \theta(v_k)\}$. In summary, for each $v_i \in V_b$, there exists an edge $e_i \in E^0$ in one of the following forms $(v_i, v_i^+)$, $(v_i^+, v_i)$, $(v_i, v_i^-)$ or $(v_i^-, v_i)$ such that
  \begin{subequations}
  \begin{align}
     \label{eqn:pf_ei1} & g\left({\begin{bmatrix} D A^T \Delta \theta\end{bmatrix}}_{e_i} \right) = 1, \;\;
      g\left({\begin{bmatrix} D A^T \Delta \phi\end{bmatrix}}_{e_i} \right) = 0 \\
     \label{eqn:pf_ei2} & e_i \neq e_j \quad \forall \; v_i \neq v_j \\
     \label{eqn:pf_ei3} & \Big| {\big\{e_i \; \vline \; v_i \in V_b \big\}} \Big| = |V_b|
  \end{align}
\end{subequations}
 Using the above argument, the inequality in (\ref{eqn:si_sibin_error_bound}) can be deduced as follows: For all feasible solutions $\Delta \theta$ of (\ref{opt:gsi}), it holds that
\begin{equation} \label{eqn:ubpf1}
\begin{array}{cl}
  & J_b - c^T g(D A^T \Delta \theta) - p^T g(A D A^T \Delta \theta) \vspace{2mm} \\
  \le & c^T g(D A^T \Delta \phi) + p^T g(A D A^T \Delta \phi) \vspace{2mm} \\
  & - c^T g(D A^T \Delta \theta) - p^T g(A D A^T \Delta \theta) 
\end{array}
\end{equation}
because $\Delta \phi$ is a feasible solution of (\ref{opt:gsi_bin}) and $J_b$ is the optimal objective value of (\ref{opt:gsi_bin}). Because of (\ref{eqn:pf_ei2}), $\big\{e_i \; \vline \; v_i \in V_b \big\}$ does not contain duplicated edges. Therefore, the right-hand-side of (\ref{eqn:ubpf1}) is equal to
\begin{equation} \label{eqn:ubpf2}
\begin{array}{l}
  \sum\limits_{v_i \in V_b} p_i \left( g\left({\begin{bmatrix} A D A^T \Delta \phi \end{bmatrix}}_{v_i} \right) -  g\left({\begin{bmatrix} A D A^T \Delta \theta \end{bmatrix}}_{v_i} \right) \right) \vspace{2mm} \\
    + \sum\limits_{\big\{e_i \; \vline \; v_i \in V_b \big\}} c_{e_i} \Bigg( g\left({\begin{bmatrix} D A^T \Delta \phi \end{bmatrix}}_{e_i} \right) - g\left({\begin{bmatrix} D A^T \Delta \theta \end{bmatrix}}_{e_i} \right) \Bigg) \vspace{2mm} \\
    + \sum\limits_{v_i \in {V^0 \setminus V_b}} p_i \left( g\left({\begin{bmatrix} A D A^T \Delta \phi \end{bmatrix}}_{v_i} \right) - g\left({\begin{bmatrix} A D A^T \Delta \theta \end{bmatrix}}_{v_i} \right) \right) \vspace{2mm} \\
    + \sum\limits_{E^0 \setminus \big\{e_i \; \vline \; v_i \in V_b \big\}} c_{e_i} \Bigg( g\left({\begin{bmatrix} D A^T \Delta \phi \end{bmatrix}}_{e_i} \right) - g\left({\begin{bmatrix} D A^T \Delta \theta \end{bmatrix}}_{e_i} \right) \Bigg).
\end{array}
\end{equation}
In addition, because of (\ref{eqn:pf_Vb_def}), (\ref{eqn:pf_ei1}) and (\ref{eqn:pf_ei3}) the expression in (\ref{eqn:ubpf2}) is equal to
\begin{equation} \label{eqn:ubpf3}
\begin{array}{l}
  \sum\limits_{v_i \in V_b} \big( p_i - c_{e_i} \big) \vspace{1mm} \\
    + \sum\limits_{v_i \in {V^0 \setminus V_b}} p_i \left( g\left({\begin{bmatrix} A D A^T \Delta \phi \end{bmatrix}}_{v_i} \right) - g\left({\begin{bmatrix} A D A^T \Delta \theta \end{bmatrix}}_{v_i} \right) \right) \vspace{2mm} \\
    + \sum\limits_{E^0 \setminus \big\{e_i \; \vline \; v_i \in V_b \big\}} c_{e_i} \Bigg( g\left({\begin{bmatrix} D A^T \Delta \phi \end{bmatrix}}_{e_i} \right) - g\left({\begin{bmatrix} D A^T \Delta \theta \end{bmatrix}}_{e_i} \right) \Bigg).
\end{array}
\end{equation}
Because of (\ref{eqn:pf_Vb_def}) and (\ref{eqn:pf_edge_cost_noincrease}), the last two sums in (\ref{eqn:ubpf3}) are nonpositive. Therefore, it holds that
\begin{equation} \label{eqn:ubpf4} 
\begin{array}{cl}
& J_b - c^T g(D A^T \Delta \theta) - p^T g(A D A^T \Delta \theta) \vspace{2mm} \\
\le & \sum\limits_{v_i \in V_b} \big( p_i - c_{e_i} \big) \vspace{1mm} \\
    \le & \sum\limits_{v_i \in V_b} \max\limits_{v_j \in N(v_i)} \big\{p_i - c_{ij} \big\} \vspace{1mm} \\
    \le & \sum\limits_{v_i \in V^0} \max\Big\{0, \max\limits_{v_j \in N(v_i)} \big\{p_i - c_{ij} \big\} \Big\}.
\end{array}
\end{equation}
Finally, since (\ref{eqn:ubpf4}) applies to all feasible solutions $\Delta \theta$ of (\ref{opt:gsi}), the upper bound in (\ref{eqn:si_sibin_error_bound}) follows.
\end{IEEEproof}
\begin{rem}
  The full measurement assumption in Proposition~\ref{thm:si_sibin_equiv} corresponds to a special case in Theorem~\ref{thm:si_sibin_equiv_general} where $c_{ij} = 2$ for all $(v_i,v_j) \in E^0$ and $p_i = 1$ for all $v_i \in V^0$. The inequalities in (\ref{eqn:si_sibin_error_bound}) imply Proposition~\ref{thm:si_sibin_equiv}.
\end{rem}
\begin{rem} \label{rem:exact_conditions}
  Theorem~\ref{thm:si_sibin_equiv_general} suggests other situations where (\ref{opt:security_index}) and (\ref{opt:security_index_binary}) are equivalent. One example is when there is a meter on each edge and there is at most one meter in each node. In this case, ${\begin{bmatrix} {P_1}^T {P_2}^T \end{bmatrix}}^T$ does not have a zero column and $P_3$ consists of subsets of rows of an identity matrix. This corresponds to $c_{ij} = 1$ and $p_i \le 1$ for all $i,j$ in Theorem~\ref{thm:si_sibin_equiv_general}. Another situation suggesting equivalence is as follows: if an edge is not metered, then its two terminal nodes are not metered either. This corresponds to a case when $p_i \le \min\limits_{v_j \in N(v_i)} c_{ij}$ for all $v_i \in V^0$, implying that $\max\limits_{v_j \in N(v_i)} \big\{p_i - c_{ij} \big\} = 0$.
\end{rem}
\begin{rem} \label{rem:approximation}
  Without the full measurement assumption or conditions such as those described in Remark~\ref{rem:exact_conditions}, solving (\ref{opt:security_index_binary}) can lead to an approximate solution to (\ref{opt:security_index}) with an error upper bound provided by (\ref{eqn:si_sibin_error_bound}). This error bound, however, is rather conservative since the summation is over all nodes $v_i \in V^0$. As developed in the proof, the summation is in fact over a subset $V_b$ of $V^0$. However, in general it is difficult to characterize the $V_b$ which leads to the tightest possible upper bound without first solving (\ref{opt:security_index}) to optimality.
\end{rem}

\subsection{Reformulating the Security Index Problem into Generalized Min Cut Problem with Costly Nodes}
\label{subsec:gmc_costly_nodes}
The above discussion suggests that the (exact or approximate) solution to the security index problem is obtained by solving (\ref{opt:gsi_bin}), whose graph interpretation will be the focus of this subsection. In (\ref{opt:gsi_bin}) the choice of 0 or 1 for each entry of $\Delta \theta$ is a partitioning of the nodes into two parts. The constraint ${A(:,e)}^T \Delta \theta \neq 0$ enforces that the two end nodes of edge $e$, denoted as $v_s$ and $v_t$, must be in two different parts of the partition. In the objective function, the term $c^T g(D A^T \Delta \theta)$ is the sum of the edge weights of the edges whose two ends are in different parts (i.e., edges that are ``cut'', in an undirected sense). In addition, since $\Delta \theta$ has binary entries, a row of $A D A^T \Delta \theta$ is zero if and only if the corresponding node and all its neighbors are in the same part of the partition (i.e., none of the incident edges are cut). Therefore, the term $p^T g(A D A^T \Delta \theta)$ in the objective function is the sum of the node weights of the nodes connected to at least one cut edge. In summary, (\ref{opt:gsi_bin}) can be reinterpreted as a generalized minimum cut problem on an undirected graph (i.e., the original power network graph with the edge direction ignored). The generalization is due to the presence of the node weights.

We now define formally the \mincutmodif problem (on any given directed graph) of which (\ref{opt:gsi_bin}) is a special case. Let $G(V,E)$ be a directed graph  (we will see that the problem can be particularized to undirected graphs), where $V$ denotes the set of nodes $\{v_1,\dots,\kcs{v_{n+1}}\}$, and $E$ the set of directed edges; and suppose that a cost $c_{ij}\geq 0$ is associated to each directed edge $(v_i,v_j)$ and a cost $p_i\geq 0$ is associated to each node $v_i$. We designate two special nodes: a source node $v_s$ and a sink node $v_t$. The problem is the following:

\begin{problem}\label{prob-mincutmodif}
\begin{equation} \label{opt:min_cost_node_partition}
\begin{array}{l}
  \textrm{\bf {The \mincutmodif problem.}} \vspace{1mm} \\
  \textrm{Find a partition of $V$, denoted as $P = \{S_s,S_t\}$, such that} \vspace{1mm} \\
  S_s,S_t \subset V, \quad S_s \cap S_t = \emptyset, \quad S_s \cup S_t = V, \quad s \in S_s, \quad t \in S_t \ \vspace{1mm} \\
  \textrm{which minimizes the cost} \vspace{1mm} \\
  \begin{array}{ccl}
  C(P) & = & \sum\limits_{(v_i,v_j) \in E: v_i \in S_s, v_j \in S_t} c_{ij} \vspace{1mm} \\ & & + \sum\limits_{v_i \in S_s : \exists (v_i,v_j) \in E: v_j \in S_t} p_i + \sum\limits_{v_j \in S_t : \exists (v_i,v_j) \in E: v_i \in S_s} p_j.
  \end{array}
\end{array}
\end{equation}
\end{problem}
By convention, if $v_i\in S_s, v_j\in S_t,$ for two nodes $v_i,v_j,$ we will say that both these nodes, and the edge $(v_i,v_j),$ are \emph{in the cut}, or that this edge is cut.

Notice that in a directed graph an edge $(v_i,v_j)$ is cut if $v_i \in S_s$ and $v_j \in S_t$ but not in the reverse case, where  $v_i \in S_t$ and $v_j \in S_s$, and the cost $c_{ij}$ is not incurred in that latter case. This asymmetry disappears however in symmetric graphs, in which to each edge $(v_i,v_j)$ with weight $c_{ij}$ corresponds a symmetric edge $(v_j,v_i)$ with same weight. For these symmetric graphs, the cost $c_{ij}$ is incurred as soon as $v_i$ and $v_j$ are not in the same set. Indeed, exactly one among $(v_i,v_j)$ and $(v_j,v_i)$ is in the cut in that case. The cost (\ref{opt:min_cost_node_partition}) consists then of the sum of the $c_{ij}$ over all pairs of nodes $v_i,v_j$ that are in different sets, and consists of the sum of the $p_i$ over all nodes that are adjacent to nodes in a different set. This is problem (\ref{opt:gsi_bin}). In addition, by letting $c_{ij}=c_{ji}=2$ for every edge and $p_i=1$ for every node, one recovers problem~(\ref{opt:security_index_binary}) under the full measurement assumption. We will show in Section~\ref{section-reformulation} how to solve Problem \ref{prob-mincutmodif}, and therefore the problems in (\ref{opt:gsi_bin}), (\ref{opt:security_index_binary}) and (\ref{opt:security_index}).



\section{An efficient solution to the tractable cases of the security index problem}
\label{section-reformulation}

This section presents the efficient solution to the \mincutmodif problem (i.e., Problem~\ref{prob-mincutmodif}) introduced in Section~\ref{subsec:gmc_costly_nodes}. The proposed solution method also solves the security index problem under the full measurement assumption, since this problem is a special case of the \mincutmodif problem.

\subsection{Construction of an Auxiliary Graph} \label{subsec:aux_graph_construction}
Consider a {directed} graph $G(V,E),\, V=\{v_1,\dots,v_{n+1}\}$ with a set of nonnegative weights $c_{ij} \ge 0$, and $p_i \ge 0$ for each node $v_i\in V$, a source node $v_s$ and a sink node $v_t$. We build an auxiliary graph $\tilde G$ using the following algorithm, illustrated in Fig.~\ref{fig:repres_graph} on an example:

\begin{figure}
\centering
\centering \scalebox{.4} {\includegraphics{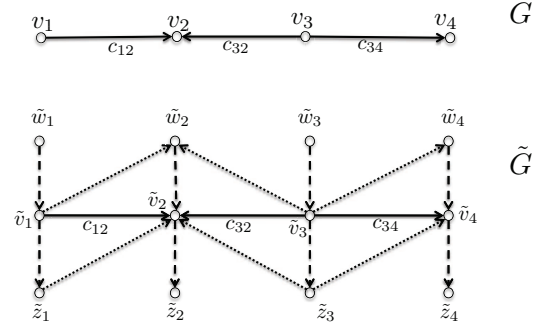}}
\caption{Representation of the auxiliary graph $\tilde G$ associated to the graph $G$. The dotted diagonal edges all have the same weight $C> \max p_i$. The vertical dashed edges linking $\tilde w_i$ to $\tilde v_i$ and $\tilde v_i$ to $\tilde z_i$ have weight $p_i$.}\label{fig:repres_graph}
\end{figure}

\begin{enumerate}
\item Define  the set $\tilde V=\{\tilde v_i,\tilde w_i, \tilde z_i: 1 \leq i\leq n+1\}$ of nodes of the auxiliary graph.
\item Designate $\tilde v_s$ and $\tilde v_t$ as source and sink nodes respectively.
\item For all $1\leq i \leq n+1$, add the two directed edges $(\tilde w_i,\tilde v_i)$ and $(\tilde v_i,\tilde z_i),$ both with cost $p_i.$
\item For all $1\leq i,j \leq n+1:\, (v_i, v_j) \in E$
  \begin{itemize} \item  add the edge $(\tilde v_i,\tilde v_j)$ with cost $c_{ij}.$
\item  add the two edges $(\tilde v_i,\tilde w_j)$ and $(\tilde z_i,\tilde v_j),$ both with a cost $C>\max_i p_i.$
\end{itemize}
\end{enumerate}

The intuition behind the construction of $\tilde G$ is the following: Suppose that one wants to cut the edge $(\tilde v_i,\tilde v_j)$, then one must also cut at least either $(\tilde v_i, \tilde z_i)$ or $(\tilde z_i, \tilde v_j)$ (see Fig.~\ref{fig:repres_graph}). Because the latter has a higher cost $C$, one will naturally cut $(\tilde v_i, \tilde z_i)$, incurring a cost $p_i$. Moreover, since that edge does not depend on $j$, one just needs to cut it (and pay the associated cost) once, independently of the number of other edges $(\tilde v_i,\tilde v_k)$ that will be cut. A similar reasoning applies to the path  $(\tilde v_i, \tilde w_j)$ or $(\tilde w_j, \tilde v_j)$. Therefore, the cost of a minimum cut on $\tilde G$ will consists of the sum of all $c_{ij}$ for all edges $(\tilde v_i, \tilde v_j)$ in the cut, and of the sum of all $p_i$ for nodes incident to one or several edges
$(\tilde v_i, \tilde v_j)$ or $(\tilde v_j, \tilde v_i)$ in the cut, i.e., to the cost of the  equivalent cut on the initial graph $G$, taking the costly nodes into account.

\subsection{Equivalence with Min Cut on the Auxiliary Graph}

We now show formally that solving the standard \mincut problem on this weighted graph provides a solution to Problem~\ref{prob-mincutmodif} on the initial graph, and that a solution is obtained by directly translating the partition of the $\tilde v_i$ into the equivalent partition of the $v_i$.

\begin{thm} \label{thm:equiv}
Consider a graph $G(V,E)$ with a set of weights $c_{ij}\geq 0$ for each edge $(v_i,v_j)\in E$, and $p_i\geq 0$ for each node $v_i\in V$, a source node $v_s$ and a sink node $v_t$. Let $\tilde G(\tilde V,\tilde E)$ be the modified graph obtained from $G$ by the procedure described above, and the partition $\tilde V = \{\tilde S_s, \tilde S_t\}$ be an optimal solution of the standard \mincut problem for $\tilde G$. Then the partition $\{S_s,S_t\}$ of $V$, obtained by letting $v_i\in S_s$ if and only if $\tilde v_i \in \tilde S_s$, is an optimal solution to Problem \ref{prob-mincutmodif} on $G$.
\end{thm}

\begin{IEEEproof}
Let us call respectively $c^*$ and $\tilde c^*$ the optimal cost of Problem \ref{prob-mincutmodif} on the graph $G$ and \mincut problem on the graph $\tilde G.$ In the sequel, we always assume that the source and sink nodes belong to the appropriate set of the partition.

We first prove that $\tilde c^* \leq  c^*,$ by showing that for any
cut in $G$ with cost $c$ (i.e., the sum of the costs of the edges {\textbf {and}} the nodes in the cut is $c$), one can build a 
cut in $\tilde G$ whose cost is equal to $c$ in the following way: For any $1\leq i\leq n+1,$
\begin{enumerate}
\item If $v_i \in S_s$, and all the out-neighbors of $v_i$ are in $S_s,$  put $\tilde v_i,$ $\tilde w_i$ and $\tilde z_i$ in $\tilde S_{s}.$
\item  if $v_i \in S_t$, and  all the in-neighbors of $v_i$ are in $S_t,$ put $\tilde v_i,$ $\tilde w_i$ and $\tilde z_i$ in $\tilde S_{ t}.$
\item if $v_i \in S_s$, and at least one  out-neighbor of $v_i$ is in $S_t,$ put $\tilde v_i,$ $\tilde w_i$ in $\tilde S_{s}$ and $\tilde z_i$ in $\tilde S_t.$
\item if $v_i \in S_t$, and at least one in-neighbor of $v_i$ is in $S_s,$ put $\tilde v_i,$ $\tilde z_i$ in $\tilde S_{ t}$ and $\tilde w_i$ in $\tilde S_{s}.$
\end{enumerate}

One can verify that no edge with cost $C$ is in the cut, and
that an edge $(\tilde v_i,\tilde v_j)$ is in the cut if and only if the corresponding edge $(v_i,v_j)$ (which has the same weight) is in the initial cut.
Moreover, for every node $i$, the edge $(\tilde w_i,\tilde v_i)$, of weight $p_i$, will be in the cut if and only if at least one edge arriving at $v_i$ was in the initial cut. Similarly, the edge $(\tilde v_i,\tilde z_i)$ will be in the cut if and only if at least one edge leaving $v_i$ is in the initial cut.
So, there will be a contribution $p_i$ to the total cost if at least an edge arriving at $v_i$ is in the cut or at least one edge leaving $v_i$ is in the cut (note that the two situations cannot happen simultaneously.)\\
As a conclusion, the cost of the 
cut $\{S_s,S_t\}$ in $G$ (counting the weights of the nodes) is equal to the cost of the 
cut $\{\tilde S_s, \tilde S_t\}$ in $\tilde G.$\\

Consider now an arbitrary 
cut in $\tilde G$, and the corresponding 
cut in $G$ obtained by putting $v_i$ in $S_s$ if and only if $\tilde v_i \in \tilde S_{s},$ as explained in the statement of this theorem. We show that the cut of $G$ obtained has a cost (taking the vertex costs $p_i$ into account) smaller than or equal to the cost of the initial 
cut. This will imply that $\tilde c^* \geq  c^*.$

The cost of this new cut $\{S_{s},S_t\}$  consists indeed of all the $c_{ij}$ of edges $(v_i,v_j)$ in the cut, and all the $p_i$ of the nodes at which arrives, or from which leaves an edge in the cut.

Consider first an edge $(v_i,v_j)$ in the cut, i.e., $v_i\in S_s,v_j\in S_t$. By construction, this implies that $\tilde v_i\in \tilde S_{s}$ and $\tilde v_j \in \tilde S_{t}$ so that the edge $(\tilde v_i, \tilde v_j)$ was also in the cut in $\tilde G$, incurring a same cost $c_{ij}$.

Consider now a node $v_i$ from which leaves at least one edge in the cut, incurring thus a cost $p_i$. (A symmetric reasoning applies if an edge in the cut arrives at $v_i$, and no node has edges in the cut both leaving from and arriving at it.) Call $v_j$ the node at which arrives that edge. We have thus $v_i\in S_s$ and $v_j\in S_t$, and therefore $\tilde v_i\in \tilde S_{s}$, $\tilde v_j \in \tilde S_{t}$ in $\tilde G$. This implies that one edge of the path consisting of $(\tilde v_i,\tilde z_i)$ and $(\tilde z_i,\tilde v_j)$ is in the cut. These edges have respective costs $p_i$ and $C>p_i$, so that a cost at least $p_i$ will be incurred by the cut in $\tilde G$. Note moreover that none of these edges will appear when considering other nodes and be counted more than once.

We have thus shown that to each cost in the cut $\{S_s,S_t\}$ for Problem \ref{prob-mincutmodif}  corresponds a larger or equal cost in $\{\tilde S_{s},\tilde S_{t}\}$ for the \mincut problem, and thus that the former has a smaller cost.

Therefore, if one takes any 
cut of optimal cost $\tilde c^*$ for the \mincut problem on $\tilde G$, and applies the procedure described in the theorem, one obtains a 
cut of $G$ with a smaller or equal cost for Problem \ref{prob-mincutmodif}. Since we have proved that the optimal cost of the latter problem is at least $\tilde c^*$, this implies that $\tilde c^* = c^*$ and that the cost obtained is optimal for Problem \ref{prob-mincutmodif} on $G$.
\end{IEEEproof}

There exist many efficient polynomial time algorithms solving the \mincut problem exactly when the weights are nonnegative (e.g. \cite{Stoer:1997:SMA:263867.263872,FF_MAX-FLOW}). Theorem \ref{thm:equiv} implies that the same algorithms can be used to solve efficiently Problem~\ref{prob-mincutmodif}, and therefore problem (\ref{opt:security_index_binary}), and problem (\ref{opt:security_index}) in the fully measured case. Moreover, observe that the size of this new graph $\tilde G$ is proportional to that of $G$, as it has $3n$ nodes and $3|E| + 2n$ edges. The order of the polynomial measuring the efficiency of the algorithms remains therefore unchanged. In particular, if the standard \mincut problem on the new graph $\tilde G$ is solved using the algorithm in \cite{Stoer:1997:SMA:263867.263872} whose complexity is $O(n|E|+n^2\log(n))$, our algorithm has the same complexity.

Finally, consider a slight generalization of Problem~\ref{prob-mincutmodif} in which each node contains two different weights (one for cutting outgoing edges and the other for cutting incoming edges). Then with a corresponding modification in the auxiliary graph construction procedure in Section~\ref{subsec:aux_graph_construction} (in the fourth bullet), the proposed method can still solve the generalization in polynomial time.

\section{The Original Security Index Problem Targeting Edge and Node} \label{sec:si_edge_only}
The relationship between the original security index problem in (\ref{opt:card_min_con}), the problem in (\ref{opt:security_index}) and its binary restriction in (\ref{opt:security_index_binary}) is summarized as follows: In the case where $H(k,:)$ in (\ref{opt:card_min_con}) corresponds to the row of $P_1 D A^T$ and $-P_2 D A^T$, (\ref{opt:card_min_con}) can be restated as (\ref{opt:security_index}) with an appropriate choice of $e$. Consequently, solving (\ref{opt:security_index_binary}) either exactly solves (\ref{opt:card_min_con}) or approximately solves (\ref{opt:card_min_con}) with an error bound provided by (\ref{eqn:si_sibin_error_bound}), depending upon whether the full measurement assumption or similar ones in Remark~\ref{rem:exact_conditions} are satisfied or not.

Next, consider the case where $H(k,:)$ corresponds to a row of $P_3 A D A^T$. The constraint $H(k,:) \dt = 1$ means that the power injection at the target node, denoted $v_s$, is nonzero. This implies that at least one edge incident to $v_s$ should have nonzero edge flow. Let $e_i$ with $i = 1,2,\ldots$ denote the column indices of $A$ of the incident edges of $v_s$. For any given $k$, consider the following instances (parameterized by $e_i$)
\begin{equation} \label{opt:card_min_con_extra_con}
  \begin{array}{cccl}
    J^i_{(\ref{opt:card_min_con_extra_con})} & \triangleq & \mathop{\textrm{min}}\limits_{\Delta \theta \in \mathbb{R}^{n+1}} & {\rm card}\big( H \Delta \theta \big) \vspace{1mm} \\
    & & \textrm{subject to} & H(k,:) \Delta \theta \neq 0 \vspace{1mm} \\
    & & & {A(:,e_i)}^T \Delta \theta \neq 0.
  \end{array}
\end{equation}
The minimum of $J^i_{(\ref{opt:card_min_con_extra_con})}$, over all $e_i$, is the optimal objective value of (\ref{opt:card_min_con}). In addition, consider a relaxation of (\ref{opt:card_min_con_extra_con}) as
\begin{equation} \label{opt:security_index_ei}
  \begin{array}{cccl}
    J^i_{(\ref{opt:security_index_ei})} & \triangleq & \mathop{\textrm{min}}\limits_{\Delta \theta \in \mathbb{R}^{n+1}} & {\rm card}\big( H \Delta \theta \big) \vspace{1mm} \\
    & & \textrm{subject to} & {A(:,e_i)}^T \Delta \theta \neq 0,
  \end{array}
\end{equation}
and its binary restriction
\begin{equation} \label{opt:security_index_binary_ei}
  \begin{array}{cccl}
    J^i_{(\ref{opt:security_index_binary_ei})} & \triangleq & \mathop{\textrm{min}}\limits_{\Delta \theta \in {\{0,1\}}^{n+1}} & {\rm card}\big( H \Delta \theta \big) \vspace{1mm} \\
    & & \textrm{subject to} & {A(:,e_i)}^T \Delta \theta \neq 0.
  \end{array}
\end{equation}
(\ref{opt:security_index_ei}) is an instance of (\ref{opt:security_index}), and the fact that (\ref{opt:security_index_ei}) has one fewer constraint than (\ref{opt:card_min_con_extra_con}) implies that
\begin{equation} \label{eqn:si_relax_lb}
  J^i_{(\ref{opt:security_index_ei})} \le J^i_{(\ref{opt:card_min_con_extra_con})}, \quad \forall \; i.
\end{equation}
On the other hand, (\ref{opt:security_index_binary_ei}) is an instance of (\ref{opt:security_index_binary}), and
\begin{equation} \label{eqn:si_bin_ub}
  J^i_{(\ref{opt:security_index_binary_ei})} \ge J^i_{(\ref{opt:card_min_con_extra_con})}, \quad \forall \; i,
\end{equation}
because if $\Delta \theta \in {\{0,1\}}^{n+1}$ is feasible to (\ref{opt:security_index_binary_ei}), then it is also feasible to (\ref{opt:card_min_con_extra_con}). Notice, however, that a feasible solution of (\ref{opt:security_index_ei}) need not be feasible to (\ref{opt:card_min_con_extra_con}). Let $i^{\star}$ be defined such that $J^{i^{\star}}_{(\ref{opt:card_min_con_extra_con})} = \min\limits_i J^i_{(\ref{opt:card_min_con_extra_con})}$. The full measurement assumption or similar ones in Remark~\ref{rem:exact_conditions} implies that $J^{i^{\star}}_{(\ref{opt:security_index_binary_ei})} = J^{i^{\star}}_{(\ref{opt:security_index_ei})}$. This, together with (\ref{eqn:si_relax_lb}) and (\ref{eqn:si_bin_ub}), suggests that
\begin{displaymath}
  J^{i^{\star}}_{(\ref{opt:card_min_con_extra_con})} \le J^{i^{\star}}_{(\ref{opt:security_index_binary_ei})} = J^{i^{\star}}_{(\ref{opt:security_index_ei})} \le J^{i^{\star}}_{(\ref{opt:card_min_con_extra_con})}.
\end{displaymath}
This implies that the equalities above hold throughout, and solving (\ref{opt:security_index_binary}) (by solving (\ref{opt:security_index_binary_ei})) indeed solves the original security index problem in (\ref{opt:card_min_con}) (by solving (\ref{opt:card_min_con_extra_con})). On the other hand, if the full measurement assumption does not hold, then
\begin{displaymath}
  J^{i^{\star}}_{(\ref{opt:card_min_con_extra_con})} \le J^{i^{\star}}_{(\ref{opt:security_index_binary_ei})} \le J^{i^{\star}}_{(\ref{opt:security_index_ei})} + \Delta J \le J^{i^{\star}}_{(\ref{opt:card_min_con_extra_con})} + \Delta J,
\end{displaymath}
where the error upper bound $\Delta J$ can be obtained from (\ref{eqn:si_sibin_error_bound}). In conclusion, all exact or approximate results pertaining to the case between (\ref{opt:security_index}) and (\ref{opt:security_index_binary}) apply to the case between (\ref{opt:card_min_con}) and (\ref{opt:security_index_binary}). As discussed in Remark~\ref{rem:approximation}, the above error bound might be conservative. The approximation quality in practice will be demonstrated in Section~\ref{section-example} containing some numerical examples on benchmark power networks.

\section{Numerical Examples} \label{section-example}

{\subsection{Simple Illustrative Example of Problem~\ref{prob-mincutmodif}} \label{subsec:simple_example}
To illustrate that the proposed method is exact while previous methods (e.g., \cite{SSJ_CDC2011_mincut,SSJ_ckt}) are not, consider an instance of Problem \ref{prob-mincutmodif} depicted in Fig.~\ref{fig:example}. Only two partitions need to be considered: $S_s = \{v_s\}$ and $S_s = \{v_s,v_1,v_2\}$, with the respective objective values being 8 and 9 (indeed, the choice $S_s = \{v_s,v_2\}$ is strictly worse than $\{v_s,v_1,v_2\}$). As a comparison, the methods in \cite{SSJ_CDC2011_mincut,SSJ_ckt} are also attempted. In particular, both \cite{SSJ_CDC2011_mincut,SSJ_ckt} solve standard \mincut problems with edge weights only. In \cite{SSJ_CDC2011_mincut} the node weights are simply ignored, while in \cite{SSJ_ckt} the node weights are indirectly accounted for by adding them to the weights of the incident edges. Table \ref{tab:cut_cost} summarizes the objective values of the source sets $\{v_s\}$ and $\{v_s,v_1,v_2\}$ for the three graph setups. The italic numbers indicate the optimal objective values in the respective methods, suggesting that both \cite{SSJ_CDC2011_mincut,SSJ_ckt} incorrectly choose $S_s = \{v_s,v_1,v_2\}$, which is suboptimal to Problem \ref{prob-mincutmodif} in the current paper.

\begin{figure}[t]
\centering
\centering \scalebox{.4} {\includegraphics{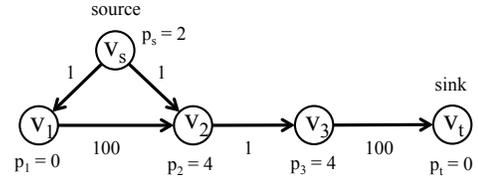}}
\caption{An instance of Problem \ref{prob-mincutmodif}. $v_s$ and $v_t$ are the source and sink nodes, respectively. The numbers next to the edges are the edge weights, and the node weights are labeled, for example, as $p_2 = 4$ for node $v_2$.}\label{fig:example}
\end{figure}

\begin{table}
  \begin{center}
    \caption{The objective values of source sets $S_s = \{v_s\}$ and $S_s = \{v_s,v_1,v_2\}$ in the graph setups of the current paper, \cite{SSJ_CDC2011_mincut} and \cite{SSJ_ckt}. As one can see, only our method finds the optimal cut.}
    \begin{tabular}{|c|c|c|c|}
      \hline
      $S_s$ & cost in our method & cost in \cite{SSJ_CDC2011_mincut} & cost in \cite{SSJ_ckt} \\
      \hline
      $\{v_s\}$ & {\emph{8}} & 2 & 10 \\
      \hline
      $\{v_s,v_1,v_2\}$ & 9 & {\emph{1}} & {\emph{9}} \\
      \hline
    \end{tabular}
    \label{tab:cut_cost}
  \end{center}
\end{table}

Constructing the auxiliary graph as described in Section~\ref{section-reformulation} and solving the corresponding standard \mincut problem leads to the node partitioning in Fig.~\ref{fig:example_aug}. In the auxiliary graph the optimal source set is $\{v_s,w_1,w_2\}$, with the objective value being 8. According to the rule in Theorem~\ref{thm:equiv}, $\{v_s\}$ is the source set returned by the proposed procedure in this paper. It correctly solves Problem \ref{prob-mincutmodif}.
\begin{figure}[t]
\centering
\centering \scalebox{.4} {\includegraphics{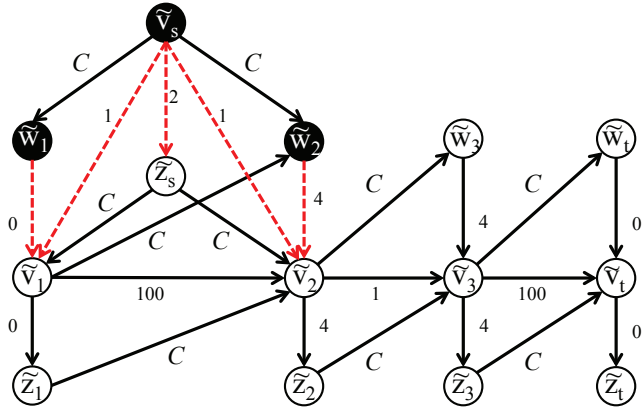}}
\caption{Solving the standard \mincut problem in the auxiliary graph corresponding to the instance in Fig.~\ref{fig:example} (the irrelevant node $w_s$ is not shown). $C$ is a large scalar constant defined in the auxiliary graph construction procedure in Section~\ref{subsec:aux_graph_construction}. The black nodes form the optimal source set (in the auxiliary graph), and the dotted red edges are cut. The optimal objective value is 8.}\label{fig:example_aug}
\end{figure}

}

\subsection{The Security Index Problem on Benchmark Systems} \label{subsec:benchmark}
\begin{figure}[tbh]
\centering
\centering \scalebox{0.34} {\includegraphics{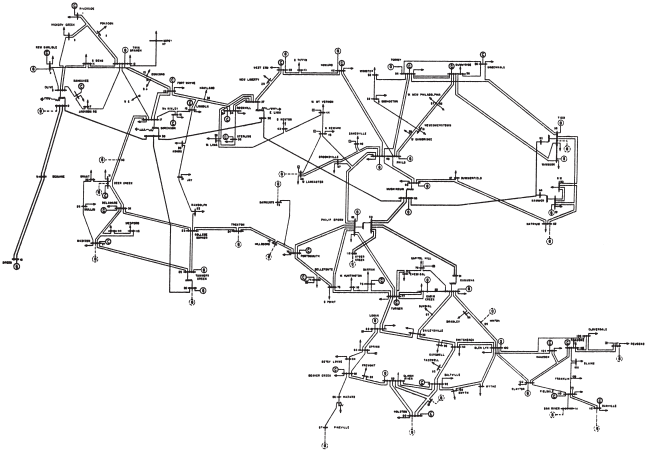}}
\caption{IEEE 118-bus benchmark system \cite{IEEE118bus}}\label{fig:118bus6}
\end{figure}
To demonstrate the effectiveness and accuracy of the proposed solution, the security index problem for two benchmark systems is considered (IEEE 118-bus \cite{IEEE118bus} and Polish 2383-bus \cite{MATPOWER}). See Fig.~\ref{fig:118bus6} for an illustration of the 118-bus system.

First, the full measurement case is considered. The security index problem in (\ref{opt:card_min_con}) is solved for each measurement, using the proposed solution and the methods from \cite{SSJ_CDC2011_mincut,SSJ_ckt}. The proposed method is guaranteed to provide the exact optimal solutions, as explained earlier in the paper. Both for the 118-bus and 2383-bus cases, the methods from \cite{SSJ_CDC2011_mincut,SSJ_ckt} are experimentally found to provide the exact solutions (though this is not guaranteed theoretically). The computation times for the three methods are listed in Table~\ref{tab:solve-time_full_meas}, indicating that the methods have similar efficiency. The guarantee of optimality provided by our approach is obtained at no additional computational cost. The computation was performed on a PC with 2.4GHz CPU and 2GB of RAM. The minimum cut problems are solved using the MATLAB Boost Graph Library \cite{Gleich06contentsmatlab,2002:BGL:504206}.
\begin{table}[!bth]
  \begin{center}
    \caption{computation times for all security indices in the full measurement case for the IEEE 118-bus and Polish 2383-bus benchmarks.}
    \begin{tabular}{|c|c|c|c|}
      \hline
      & our method & \cite{SSJ_CDC2011_mincut} & \cite{SSJ_ckt} \\
      \hline
      IEEE 118-bus & 0.18s & 0.23s & 0.24s \\
      \hline
      2383-bus & 29s & 29s & 33s \\
      \hline
    \end{tabular}
    \label{tab:solve-time_full_meas}
  \end{center}
\end{table}

\begin{figure}[t]
\centering
\centering \scalebox{0.6} {\includegraphics{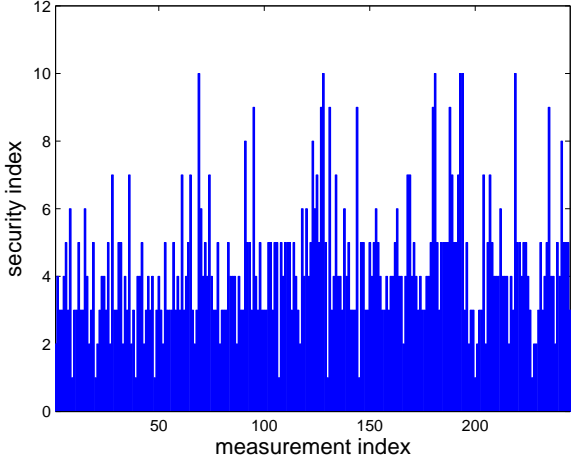}}
\caption{Security indices for the partially measured 118-bus system. Security indices are computed using the big $M$ method with $M = 10^4$.}\label{fig:si118}
\end{figure}

Next, (\ref{opt:card_min_con}) is considered when the full measurement assumption is removed. That is, the matrices $P_1$, $P_2$ and $P_3$ in (\ref{def:H_matrix}) need not be identities. In this test, the 118-bus system is considered. In the measurement system about 50\% of power injections and power flows are measured. The measurements are chosen randomly, and the measurement system is verified to be observable (i.e., the corresponding $H_{2:}$ has full column rank ($=n$)). Since (\ref{opt:card_min_con}) is NP-hard in general,  no efficient solution algorithm has been known. Enumerative algorithms include, for instance, enumeration on the support of $H \dt$, finding the maximum feasible subsystem for an appropriately constructed system of infeasible inequalities \cite{Jokar:2008:EAS:1461600.1461607}, and the big $M$ method to be described. The authors' implementations of the first two methods turn out to be too inefficient for the applications concerned. Therefore, the big $M$ method is used, which sets up and solves the following optimization problem:
\begin{equation} \label{opt:si_MILP}
  \begin{array}{cccl}
    \mathop{\rm minimize}\limits_{\Delta \theta, \; y} & \quad \sum\limits_j y(j) & & \\
    \textrm{subject to} & H \Delta \theta & \le & M y \\
    & -H \Delta \theta & \le & M y \\
    & H(k,:) \Delta \theta & = & 1 \\
    & y(j) & \in & \{0,1\} \quad \forall \; j.
  \end{array}
\end{equation}
In (\ref{opt:si_MILP}), $M$ is a user-defined constant. If $M \ge {\|H \Delta \theta^\star\|}_\infty$ for at least one optimal solution $\Delta \theta^\star$ of (\ref{opt:card_min_con}), then (\ref{opt:si_MILP}) provides the exact solution to (\ref{opt:card_min_con}). Otherwise, solving (\ref{opt:si_MILP}) yields a suboptimal solution, optimal among all solutions $\Delta \theta$ such that ${\|H \Delta \theta\|}_\infty \le M$. In principle a sufficiently large $M$ can be found to ensure that the big $M$ method indeed provides the optimal solution to (\ref{opt:card_min_con}) \cite{Schrijver:1986:TLI:17634}. However, this choice of $M$ is typically too large to be practical. In the numerical example in this section, $M$ is simply chosen to be $10^4$. (\ref{opt:si_MILP}) can be solved as a mixed integer linear program \cite{BT97} using solvers such as CPLEX \cite{CPLEX}. The solutions by the big $M$ method are treated as references for accuracy for the rest of the case study. Fig.~\ref{fig:si118} shows the (big $M$) security indices for all chosen measurements.

\begin{figure}[!t]
\centering
\centering \scalebox{0.6} {\includegraphics{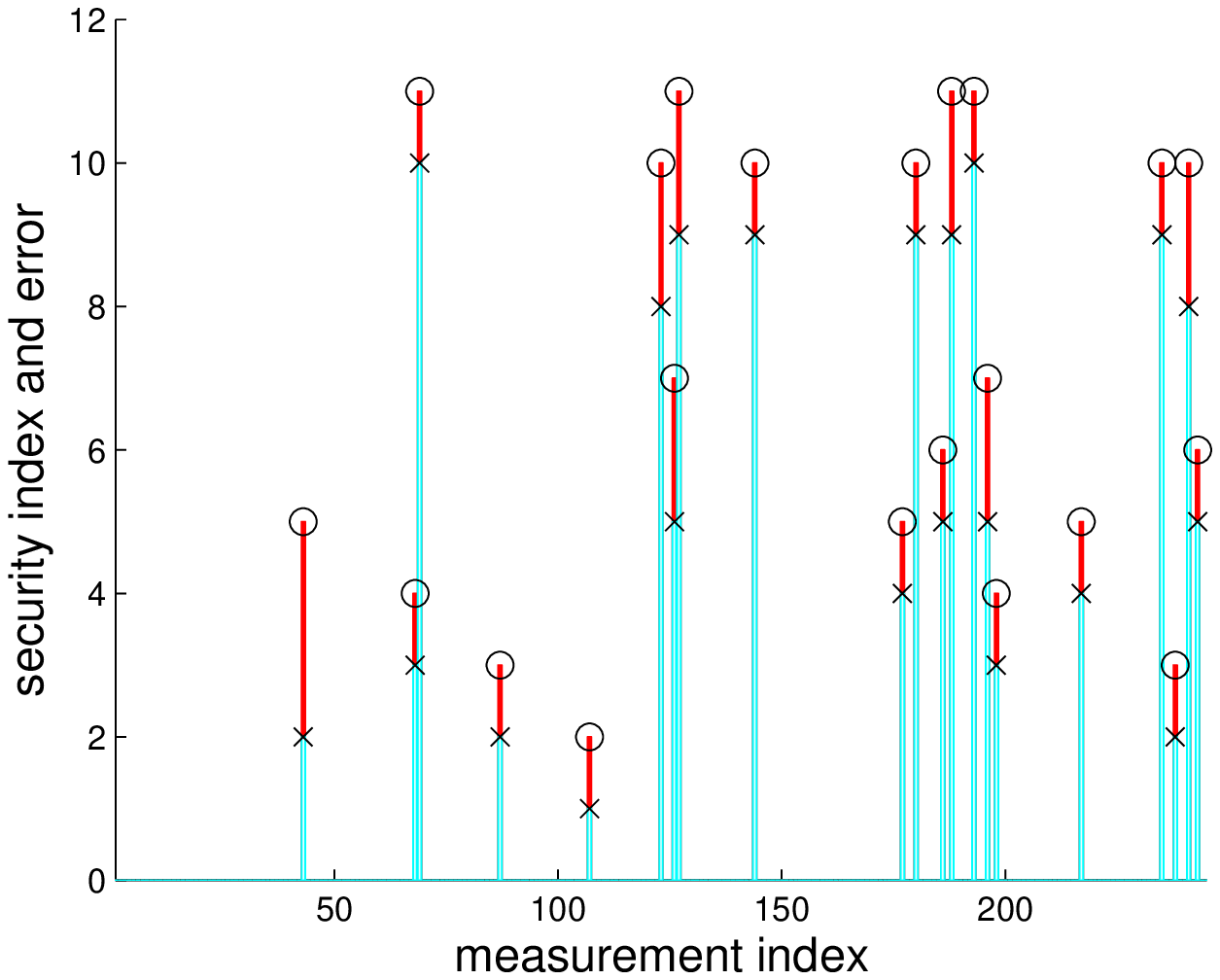}}
\caption{Security index estimates for the partially measured 118-bus system. Security index estimates are computed using the method in \cite{SSJ_CDC2011_mincut}. The figure shows only the inexact security index estimates (in circles) and the corresponding ones by the big $M$ method (in crosses). }\label{fig:si118_mc}
\end{figure}

\begin{figure}[t]
\centering
\centering \scalebox{0.6} {\includegraphics{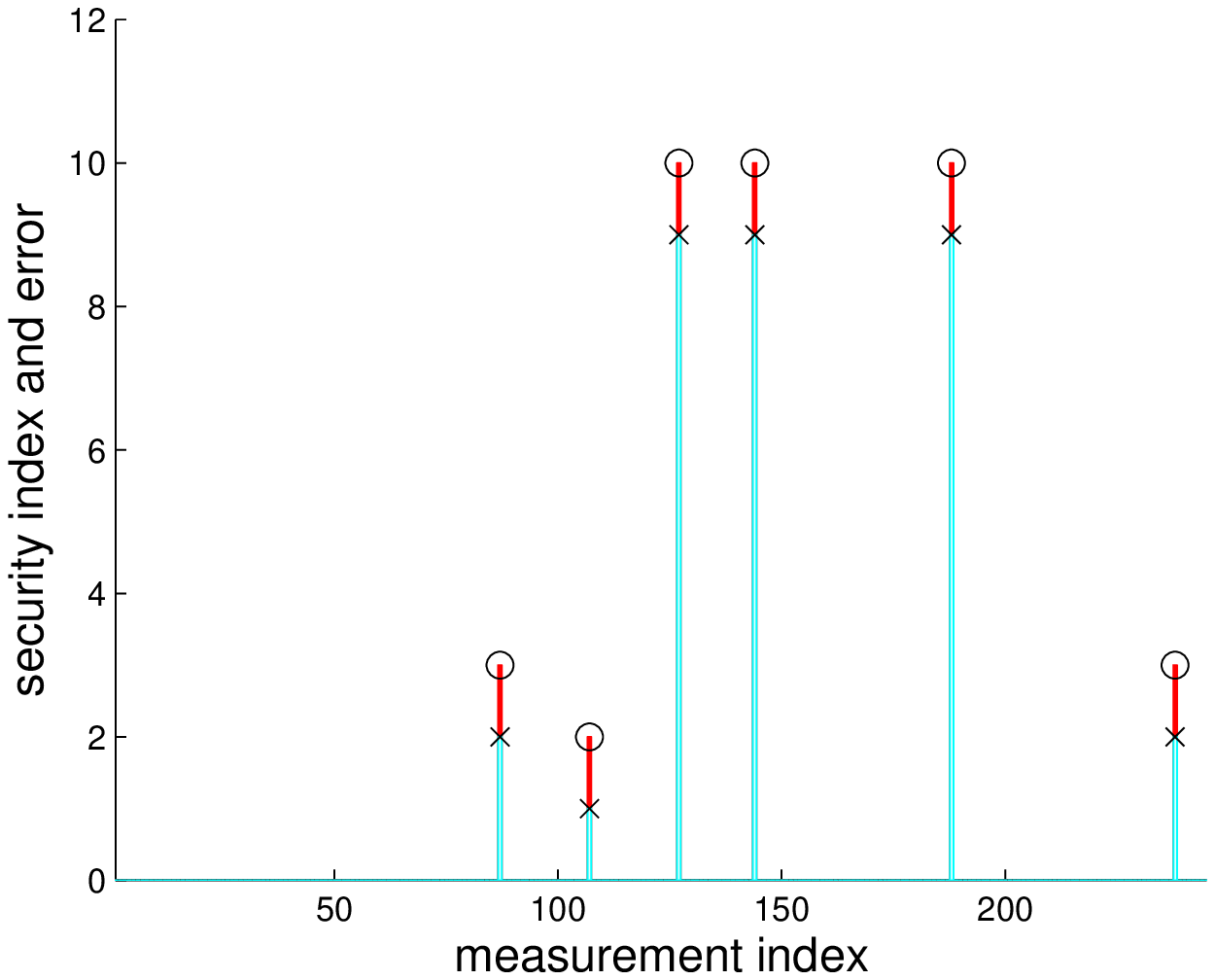}}
\caption{Security index estimates for the partially measured 118-bus system. Security index estimates are computed using the method in \cite{SSJ_ckt}. The figure shows only the inexact security index estimates (in circles) and the corresponding ones by the big $M$ method (in crosses).}\label{fig:si118_mc1}
\end{figure}

\begin{figure}[!tbh]
\centering
\centering \scalebox{0.6} {\includegraphics{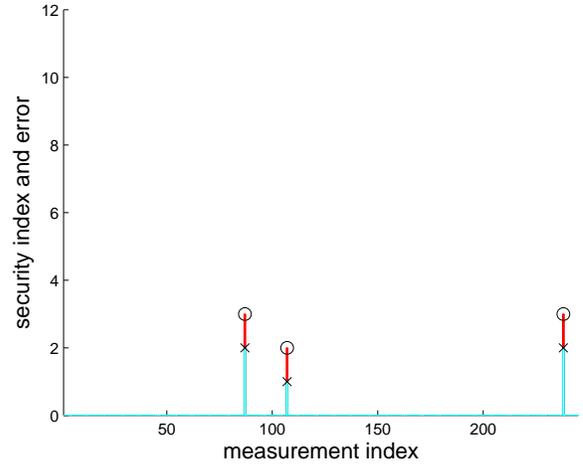}}
\caption{Security index estimates for the partially measured 118-bus system. Security index estimates are computed using the method proposed in this paper. The figure shows only the inexact security index estimates (in circles) and the corresponding ones by the big $M$ method (in crosses).}\label{fig:si118_gmc}
\end{figure}

Alternatively, as described in Section~\ref{sec:si_edge_only} a suboptimal solution to (\ref{opt:card_min_con}) can be obtained by solving (\ref{opt:security_index_binary}) exactly using the proposed method in Section~\ref{section-reformulation} or the ones from \cite{SSJ_CDC2011_mincut,SSJ_ckt}. As explained earlier, (\ref{opt:security_index_binary}) can be formulated as Problem~\ref{prob-mincutmodif} with $p_i = 1$ if and only if the injection measurement at bus $v_i$ is taken and $c_{ij} = c_{ji} \in \{0,1,2\}$ being the total number of line power flow meters on the line connecting buses $v_i$ and $v_j$. Fig.~\ref{fig:si118_mc}, Fig.~\ref{fig:si118_mc1} and Fig.~\ref{fig:si118_gmc}, respectively, show the security index test results with the three \mincut based methods (i.e., \cite{SSJ_CDC2011_mincut,SSJ_ckt} and the one proposed in this paper).
These figures show only the big $M$ security indices (in light blue, the heights of the crosses) and the overestimation (in red, the heights between the crosses and the circles) for the measurements where the \mincut based methods do not agree with big $M$. The case study indicates that, among the three \mincut based methods, the proposed method provides the most accurate suboptimal solutions to (\ref{opt:card_min_con}). In terms of computation time, the proposed method is most efficient as suggested by Table~\ref{tab:solve-time_partial_meas}.

\begin{table}[th]
  \begin{center}
    \caption{computation times for all security indices in the partial measurement case for the IEEE 118-bus benchmark.}
    \begin{tabular}{|c|c|c|c|c|}
      \hline
      & our method (s) & \cite{SSJ_CDC2011_mincut} & \cite{SSJ_ckt} & big $M$ \\
      \hline
      IEEE 118-bus & 0.17s & 4.4s & 0.21s & 118s \\
      \hline
    \end{tabular}
    \label{tab:solve-time_partial_meas}
  \end{center}
\end{table}

\section{Conclusions}\label{section-concl}
It has been assumed that the security index problem, formulated as a cardinality minimization problem, cannot be solved efficiently. This paper formally confirms this conjecture by showing that the security index problem is indeed NP-hard. Nevertheless, the security index problem can be shown to be reducible to a \mincutmodif problem (Problem \ref{prob-mincutmodif}) under the full measurement assumption. In this paper, we show that this problem is equivalent to a standard \mincut problem on an auxiliary graph of proportional size, and can therefore be solved exactly and efficiently using standard techniques for the \mincut problem. Under the full measurement assumption, this allows computing the minimal number of measurements with which one must tamper in order to feed incorrect information on the SCADA system without being detected by a BDD method. The knowledge of this number can help strategically assigning protection resources (e.g., \cite{DS_SGC2010,VSDS_jsac11}). Our method also solves a mathematically equivalent problem of robustness of the observability properties of the system with respect to the failure of some measurements, assuming again full measurement. It remains to be determined if the solution could be efficiently approximated in the general (not fully measured) case. Indeed, even though our approach already provides an approximate solution to such general problems we do not know if this approximation comes with any guarantee of accuracy.

Another interesting issue is the design question: in view of the exact solution of the security index problem presented in this paper, could one build efficient design methods in order to optimize the security index under some natural constraints?



\bibliographystyle{IEEEtran}
\bibliography{gmc_tsg,exact_l1,references_all}

\end{document}